\documentclass[11pt, reqno]{amsart}
\usepackage{amsmath, amsthm, amscd, amsfonts, amssymb, graphicx, color}
\usepackage[bookmarksnumbered,plainpages]{hyperref}

\newcommand{\bea}{\begin{eqnarray*}}
\newcommand{\eea}{\end{eqnarray*}}
\newcommand{\A}{\mathfrak A}
\newtheorem{theorem}{Theorem}[section]
\newtheorem{lemma}[theorem]{Lemma}
\newtheorem{proposition}[theorem]{Proposition}
\newtheorem{corollary}[theorem]{Corollary}
\theoremstyle{definition}

\newtheorem{example}[theorem]{Example}

\theoremstyle{remark}
\newtheorem{remark}[theorem]{Remark}
\numberwithin{equation}{section}

\begin{document}
\title{ Weighted semigroup algebras as dual Banach algebras}{}
\author{M. Abolghasemi}
\address{Department of Mathematics, Isfahan University, Isfahan, Iran}\email{m.abolghasemi@sci.ui.ac.ir}
\author{A. Rejali}
\address{Department of Mathematics, Isfahan University, Isfahan, Iran.}\email{Rejali@sci.ui.ac.ir}
\author{H.\ R.\ E.\ Vishki}
\address{Department of Pure Mathematics and Centre of Excellence
in Analysis on Algebraic Structures (CEAAS), Ferdowsi University of Mashhad\\
P. O. Box 1159, Mashhad 91775, Iran.} \email{vishki@ferdowsi.um.ac.ir}

\subjclass[2000]{43A10, 43A20, 46H20, 46H25}
\keywords{Weakly cancellative
semigroup, weighted semigroup algebra, measure algebra, dual Banach algebras. }
\begin{abstract} In this paper, among other things,   we study those   conditions under which the  weighted semigroup algebra $\ell^1(S,\omega)$ is a dual Banach algebra with respect to  predual $c_0(S)$.  Some useful examples, illustrating the results, are also included.
\end{abstract}
\maketitle

\section{Introduction and Preliminary results}

For a Banach algebra $\mathfrak{A}$ assume that ${\mathfrak{A}}^*$ and
$\mathfrak{A}^{**}$ are the first and second duals of
$\mathfrak{A}$, respectively.
For $a, b\in \mathfrak{A}$ and $f\in \mathfrak{A}^*$, we define the elements $f\cdot a$ and $a\cdot f$ of $\mathfrak{A}^*$
by the formulae
$$\langle  f\cdot a, b\rangle=\langle f ,ab\rangle~~\emph{\emph{and}}~~\langle a\cdot f, b\rangle=\langle f, ba\rangle.$$
These are the natural module actions of $\mathfrak{A}$ on $\mathfrak{A}^*$ which turns $\mathfrak{A}^*$ into a Banach $\mathfrak{A}-$module.
Also we define the so-called Arens products  $\square$ and $\lozenge $ on $\mathfrak{A}^{**}$ by the formulae:
$$\langle m~ \square ~n, f\rangle=\langle m, n\cdot f\rangle~~\emph{\emph{and}}~~\langle m~ \lozenge
~n, f\rangle=\langle n, f\cdot m\rangle, ~~~~~~~~~~~~~~~~~~~~~~~~~~~~~~~~~~~(m, n\in {\mathcal A}^{**});$$
in which,
$$\langle f\cdot m, a\rangle=\langle m, a\cdot f\rangle~~\emph{\emph{and}}~~\langle m\cdot f, a\rangle=\langle m, f\cdot a\rangle.$$
The Banach algebra $\mathfrak{A}$ is said to be Arens regular if  $``~\Box=\lozenge ~"$ on the whole of $\mathfrak{A}^{**}$. This is the case, for example, for all $C^*-$algebras, \cite {CY}, and also for $L^1(G)$ if (and only if) $G$ is finite, \cite Y. The interested reader may refer to \cite {DL} for ample information about the   Arens regularity problem.\\

A Banach algebra $\mathfrak{A}$ which is a dual space, with predual  $\mathfrak{A}_*$, is called a left (respectively, right) dual Banach algebra  if  $\mathfrak{A}_*$ is a closed left (respectively, right) $\mathfrak{A}-$submodule  of $\mathfrak{A}^*$. If $\A$ is both left and right dual Banach algebra, then it is called a dual Banach algebra. In each case we involve to the weak$^*-$topology $\sigma(\mathfrak{A}, \mathfrak{A}_*)$ on $\mathfrak{A}.$ It is easy to verify that  $\mathfrak{A}$ is a  dual Banach algebra if and only if it is a dual Banach space whose  multiplication is separately $w^*-$continuous (see, {\cite[Lemma 5.1]{DPV}}).

 It should be mentioned that the  predual of a dual Banach algebra need not be unique, in general; for instance, one may furnish a Banach space,  which admits   two distinct preduals (such as,  $\ell^1$ which  admits $c$ and $c_0$ as two distinct preduals) with the  zero product to obtain a dual Banach algebra with two distinct preduals. Therefore to avoid the risk of ambiguity we  usually determine the involved predual.
 We wish to thank M. Daws for bringing the preprint~\cite{DPW} to our attention; in which the authors studied the uniqueness of predual for a wide variety of group algebras.

We commence with the next result.
\begin{proposition}
Let $\mathfrak A$ be a  dual Banach algebra. If $\mathfrak A$ has a  bounded approximate identity then it has an identity.
\end{proposition}
\begin{proof}
 Let $(e_\alpha)$ be a  bounded approximate identity for
$\mathfrak{A}$. Clearly it has a $w^*-$cluster point $e\in ({{\mathfrak{A}}_*})^*$. By passing to a subnet, if necessary, assume that $w^*-\lim e_\alpha =e$. Consequently, for each $a\in\A$,  $w^*-\lim e_\alpha a=ea$; on the other hand $w^*-\lim e_\alpha a=a$. Therefore
$ea=a$ and so e is a left identity of $\mathfrak{A}$. A similar argument may apply for the right identity.
\end{proof}
It is well known that every von Newmann algebra is a dual Banach algebra; more precisely, a $C^*-$algebra is a dual Banach algebra if and only if it is a von Neumann algebra (see, {\cite[Example 4.4.2(c)]{Runde1}}). Also the measure algebra $ M(G)$ of a locally compact group $G$ is a dual Banach algebra (with predual $C_0(G)$). In {\cite[Theorem 4.6]{Dales-Lau-s}} it has been shown that the semigroup algebra $\ell^1(S)$ is a dual Banach algebra with predual $c_0(S)$ if and if $S$ is weakly cancellative.

 A direct verification reveals that $\mathfrak{A}^*$ is a left  $(\mathfrak{A}^{**}, \square)-$submodule (respectively, right $(\mathfrak{A}^{**}, \lozenge)-$ submodule) of $\mathfrak{A}^{***}$; therefore   $(\mathfrak{A}^{**}, \square)$ is a left dual Banach algebra as well as  $(\mathfrak{A}^{**}, \lozenge)$ is a right  dual Banach algebra (with predual $\mathfrak{A}^*$). Moreover the following result (a version of which is given in {\cite[Proposition 5.2]{DPV}}; and also a  module action version of which is given in {\cite[Proposition 3.1]{MV}}) determines the situation under which each of $(\mathfrak{A}^{**}, \square)$ and  $(\mathfrak{A}^{**}, \lozenge)$ is actually a dual Banach algebra.
\begin{proposition} [See {\cite[Proposition 5.2]{DPV}}] For a Banach algebra $\mathfrak{A}$ the following statements are equivalent.

$(i)$ $(\mathfrak{A}^{**}, \square)$  is a dual Banach algebra (with predual $\mathfrak{A}^*$).

$(ii)$  $(\mathfrak{A}^{**}, \lozenge)$ is a dual Banach algebra (with predual $\mathfrak{A}^*$).

$(iii)$ $\mathfrak{A}$ is Arens regular.
\end{proposition}
\noindent Furthermore, as it is  emphasized in {\cite[Theorem 5.1]{DPW}}, in the case where $\mathfrak A$ is Arens regular then  $\mathfrak A^*$ is the unique predual of the dual Banach algebra $\mathfrak A^{**}$, provided that  $\mathfrak A$ is an ideal in $\mathfrak A^{**}$ and $\mathfrak A^{**}$ is unital.\\

The central goal of this paper is to provide those conditions under which the  weighted
semigroup algebra $\ell^1(S, \omega)$ and semigroup measure algebra
$M_b(S)$ are dual Banach algebras.  We show that   $\ell^1(S,\omega)$ is a dual
Banach algebra with respect to the canonical predual $c_0(S)$ if and only if $\frac{1_{st^{-1}}}{\omega}\in
c_0(S)$ and $\frac{1_{t^{-1}s}}{\omega}\in c_0(S),$ for all
$s,t\in S$. From this some older results concerning the dual semigroup algebra $\ell^1(S)$ are improved.  For a locally compact Hausdorff  topological semigroup $S$, the measure algebra $M_b(S)$---from the dual Banach algebra point of view---is also studied .
\\

\section{ Weighted semigroup algebra as a dual Banach algebra }
Let $F$ and $K$  be  nonempty subsets of a semigroup  $S$ and $s\in
S$. We put

\[s^{-1}F=\{t\in S:~st\in F\},~
\emph{\emph{and}} ~Fs^{-1}=\{t\in S:~ts\in F\},\] and we also  write
$s^{-1}t$ for the set $s^{-1}\{t\}$, $FK^{-1}$ for $\bigcup\{Fs^{-1}: s\in K\}$ and $K^{-1}F$ for $\bigcup\{s^{-1}F: s\in K\}$.\\
 A semigroup $S$ is said to be  left (respectively,  right) weakly cancellative if $s^{-1}F$ (respectively, $Fs^{-1}$) is finite for each $s\in S$ and each finite subset $F$ of $S$. A semigroup $S$ is said to be weakly cancellative if it is both  left and  right weakly  cancellative.

A weight on $S$ is a function $\omega: S\rightarrow (0, \infty)$ such that
\[\omega(st)\leq \omega(s)\omega(t), ~~~~~~~(s, t\in S).\]

Let  $\omega$ be a weight on $S$, then the space
$$\ell^1(S,\omega)=\{f:~f=\sum_{s\in S} f(s)\delta_s,~~\|f\|_{1,\omega}=\sum_{s\in S} |f(s)|\omega(s)<\infty\},$$
(where, $\delta_s\in \ell^1(S,w)$ is the point mass at $s$ which can be thought as the  vector basis element of $\ell^1(S,w)$) equipped with the multiplication
\begin{eqnarray*}
f*g =(\sum_{s\in S} f(s)\delta_s)*(\sum_{t\in S}
g(t)\delta_t)=\sum_{r\in S}(\sum_{st=r}
f(s)g(t))\delta_r~,
\end{eqnarray*}
(and also define $f*g=0$, if for each $r\in S$ the equation  $st=r$ has no solution;)
is a Banach algebra which will be  called a weighted semigroup algebra.\\

The dual space of $\ell^1(S,\omega)$ is the space
$$\ell^{\infty}(S,\omega)=\{\phi:~\phi=\sum_{s\in S} \phi(s) e_s,~~\|\phi\|_{\infty,\omega}=\sup_{s\in S}\mid \frac{\phi(s)}{\omega(s)}\mid<\infty\};$$
(here we use the same point mass $e_s$ as the  vector basis element of   $l^{\infty}(S,w)$), whose  duality is given by
\[\langle \phi, f\rangle=\sum_{s\in S}\phi(s)f(s)\ \ \ (\phi\in\ell^\infty(S,\omega), f\in\ell^1(S,\omega)).\] In particular, $\langle e_s, \delta_t\rangle=\delta_{s,t}$,
where $\delta_{s,t}$ is the Kronecker delta, that is,
\begin{eqnarray*}\delta_{s,t}=&\left\{%
\begin{array}{ll}
    1 & \hbox{$s=t$} \\
    0 & \hbox{$s\neq t$}. \\
\end{array}%
\right.
\end{eqnarray*}
 From now on we write $1_F$ for the characteristic function of $F$.\\

 The  next lemma describes explicitly the module actions of $\mathfrak{A}=\ell^1(S,\omega)$ on  $\mathfrak{A^*}=\ell^{\infty}(S,\omega).$
\begin{lemma}\label{3.2}
Let $S$ be a semigroup, then
$e_s\cdot\delta_t=1_{t^{-1}s}~~ and~~~ \delta_t\cdot e_s=1_{st^{-1}},$  for every $s,t\in S.$
\end{lemma}
\begin{proof}
Let $r\in S$, then
\bea
\langle e_s\cdot\delta_t, \delta_r\rangle &=& \langle e_s,
\delta_t*\delta_r\rangle\\&=&\langle e_s, \delta_{tr}\rangle\\ &=&\delta_{s, tr}
\\&=&\{\begin{array}{c}
  1~~~~r\in t^{-1}s \\
  0~~~~r\notin t^{-1}s \\
\end{array}\\&=&\langle 1_{t^{-1}s}, \delta_r\rangle.
\eea
Therefore $e_s\cdot\delta_t=1_{t^{-1}s},$ and similarly
$\delta_t\cdot e_s=1_{st^{-1}}.$
\end{proof}
By the latter lemma, for  $f\in\ell^1(S,\omega)$ and $\phi\in\ell^{\infty}(S,\omega)$, with  $f=\sum_{s\in S} f(s)\delta_s$ and $\phi=\sum_{s\in S} \phi(s) e_s$, we have:
\[ \phi\cdot f=\sum_{s,t\in S}\phi(s)f(t)1_{t^{-1}s}, \ \ \  {\rm and}\ \  f\cdot\phi=\sum_{s,t\in S}\phi(s)f(t)1_{st^{-1}}.\]\\

We also define
$$c_0(S,\omega)=\{\phi:~\phi\in \ell^{\infty}(S,\omega) ~\emph{\emph{and}}~
\frac{\phi}{\omega}\in c_0(S)\}.$$
Then $c_0(S,\omega)$ is a closed subspace of $\ell^{\infty}(S,\omega)$; indeed, $c_0(S,\omega)$ is the closed linear span of $\{e_s : s\in S\}$ in $\ell^{\infty}(S,\omega)$. The dual space of $c_0(S,\omega)$ is $\ell^1(S,\omega)$ (whose duality is the restriction of that we have presented between  $\ell^{\infty}(S,\omega)$ and $\ell^1(S,\omega)$), and so its second dual can be identified with $\ell^{\infty}(S,\omega)$.\\

The next result  determines when  a weighted semigroup algebra is a dual Banach algebra.
\begin{theorem}\label{p2} A weighted semigroup algebra    $\ell^1(S,\omega)$ is a left (respectively, right) dual Banach algebra with predual $c_0(S,\omega)$ if and only if  $\frac{1_{st^{-1}}}{\omega}\in
c_0(S)$  (respectively, $\frac{1_{t^{-1}s}}{\omega}\in c_0(S)$), for all $s,t\in S.$
\end{theorem}
\begin{proof} Let  $\ell^1(S,\omega)$ be a right  dual Banach algebra
with predual $c_0(S,\omega)$, that is,
$$c_0(S,\omega)\cdot \ell^1(S,\omega) \subseteq c_0(S,\omega).$$
For $s,t\in S$ since $e_s\in c_0(S,\omega)$ and $\delta_t \in \ell^1(S,\omega)$, the latter inclusion together with the  equality $1_{t^{-1}s}=e_s\cdot\delta_t$,
Lemma \ref{3.2}, imply that $1_{t^{-1}s}\in c_0(S, \omega)$.

For the converse,   suppose that $f\in \ell^1(S,\omega)$ and
$\phi\in c_0(S,\omega)$. For each  $\epsilon>0$ there exist functions $\phi_0$ and $f_0$ with finite supports such that\[\|\phi-\phi_0\|_{\infty,\omega}<\epsilon ~~~~{\rm{and}}~~~~\|f-f_0\|_{1,\omega}<\epsilon.\]
Now  from the inequality  $\|\phi\cdot f\|_{\infty,\omega}\leq \|\phi\|_{\infty,\omega}\|f\|_{1,\omega};$ we have,
\begin{eqnarray*}
|\phi\cdot f| &\leq &\epsilon(\|\phi_0\|_{\infty,\omega}+\|f\|_{1,\omega}) + |\phi_0\cdot f_0|.
\end{eqnarray*}
As $\phi_0\cdot f_0=\sum_{s\in S, t\in
S}\phi_0(s)f_0(t)1_{t^{-1}s}$ is a finite combination of the functions  $1_{t^{-1}s}$ which are in  $ c_0(S,\omega)$, we have $\phi_0\cdot f_0\in c_0(S, w)$ and the latter inequality implies that $\phi\cdot f\in c_0(S, w)$, as required. The left case needs a similar proof. \\
 \end{proof}
\begin{remark}
Let $\mathfrak A$ be a Banach algebra with two preduals  ${\mathfrak A_*}_1$ and ${\mathfrak A_*}_2$. It may be that  $\mathfrak A$ is a dual Banach algebra with respect to ${\mathfrak A_*}_1$, but is not with respect to  ${\mathfrak A_*}_2$; (see,  {\cite[Example 7.30]{Dales-Lau-s}} or Example~\ref{4} below).  If ${\mathfrak A_*}_1$ and  ${\mathfrak A_*}_2$ are  isometrically isomorphic  (as Banach spaces) then a routine argument shows that $\mathfrak A$ is dual with respect to ${\mathfrak A_*}_1$ if and only if it is dual with respect to ${\mathfrak A_*}_2$. We mention that the $\mathfrak A-$module actions $\cdot$ and $\bullet$ of   ${\mathfrak A_*}_1$ and  ${\mathfrak A_*}_2$, respectively,  are related to each other via the equations
\[a\bullet x=\theta(a\cdot \theta^{-1}(x))\ \ {\rm and}\ \ x\bullet a=\theta(\theta^{-1}(x)\cdot a), \ \ \ (a\in \mathfrak A, x\in{\mathfrak A_*}_2);\]
where, $\theta :{\mathfrak A_*}_1\rightarrow{\mathfrak A_*}_2$ is the involved isometric isomorphism.

With this observations in mind, as for every two weights $\omega$ and $\omega^{'}$ on $S$,  the Banach spaces $c_0(S, \omega^{'})$ and $c_0(S, \omega)$ are isometrically isomorphic (under the mapping $\phi\mapsto \phi\frac{\omega}{\omega^{'}}$);
 $\ell^1(S,\omega)$ is a dual Banach algebra with respect to the predual $c_0(S, \omega)$ if and only if it is dual with respect to the predual $c_0(S, \omega^{'})$; (or  equivalently,  with respect to the predual $c_0(S)$). Based on these facts,  being $\ell^1(S,\omega)$ as a  dual Banach algebra
with predual $c_0(S,\omega^{'})$  depends  merely on the behavior of the weight $\omega$ imposed on  $\ell^1(S,\omega)$ and  is  quite independent from the behavior of the weight $\omega^{'}$ imposed on the predual. So for the sake of simplicity, without loss of generality we only need  to study $\ell^1(S,\omega)$ from the dual Banach algebra point of view with respect to the canonical  predual $c_0(S)$.
 \end{remark}

Based on what we have remarked in the preceding  remark, as an   immediate consequence of the latter theorem we have,
\begin{corollary}\label{p3} A weighted semigroup algebra $\ell^1(S,\omega)$ is a  dual Banach algebra with predual $c_0(S)$ if and only if $\frac{1_{st^{-1}}}{\omega}, \frac{1_{t^{-1}s}}{\omega}\in c_0(S)$, for all $s,t\in S.$
\end{corollary}

We apply  Corollary~\ref{p3} to obtain the next result, whose sufficiency  can also be found in {\cite[Proposition 5.1]{Daws}}, with a proof slightly different from ours.
\begin{proposition}\label{pp3} If $S$ is weakly cancellative then
 $\ell^1(S, \omega)$ is a  dual Banach algebra with predual $c_0(S)$. The converse is also valid in the case where $\omega$ is bounded.
\end{proposition}
\begin{proof}
 Let $S$ be weakly cancellative, then for each $s, t\in S$ the functions   $\frac{1_{t^{-1}s}}{w}$ and $\frac{1_{st^{-1}}}{w}$ vanish on the complement of $t^{-1}s\cap st^{-1}$, and the latter is finite  by the weak cancellativity. Now Corollary~\ref{p3} implies that $\ell^1(S, \omega)$ is a  dual Banach algebra with predual $c_0(S)$.
For the converse assume that  $\omega$ is bounded (by $M>0$). As $\ell^1(S, \omega)$ is a  dual Banach algebra with predual $c_0(S)$, again by Corollary~\ref{p3} there exists a
finite subset $F$ of $S$ such that for each $x\notin F$,
$$\frac{1_{t^{-1}s}}{\omega}(x)<\frac{1}{2M}.$$ Therefore $x\notin
t^{-1}s$, that is, $t^{-1}s\subseteq F,$ hence $S$ is left weakly cancellative. A similar argument may apply for the right cancellativity of $S$.
\end{proof}
Using the latter result for the special case $\omega=1,$ we present the next result of Dales {\it et al.}, \cite{Dales-Lau-s}.
\begin{corollary}[{\cite[Theorem 4.6]{Dales-Lau-s}}] For a semigroup $S$, $\ell^1(S)$ is a dual Banach algebra with the predual $c_0(S)$ if and only if $S$ is weakly cancellative.
\end{corollary}

It should be mentioned that however the weak cancellativity of $S$, by Proposition~\ref{pp3}, ensures that $\ell^1(S, \omega)$ is a dual Banach algebra with respect to all preduals $c_0(S)$, but the same result may hold in non-weakly cancellative setting, provided that the weight $\omega$  enjoys some additional conditions. The next result provides some conditions  under which  $\ell^1(S, \omega)$ is a dual Banach algebra with predual $c_0(S)$.
\begin{proposition}\label{p4} In either of the following cases $\ell^1(S, \omega)$ is a  dual Banach algebra with predual $c_0(S)$.

$(i)$ $w_l, w_r\in c_0(S)$; where,  $w_l(t)=\sup_{s\in S}|\frac{\omega(st)}{\omega(t)}|$ and $w_r(t)=\sup_{s\in S}|\frac{\omega(ts)}{\omega(t)}|$, $(t\in S).$

$(ii)$ $\frac{1}{\omega}\in c_0(S).$
\end{proposition}
\begin{proof} $(i)$ Assume that  $w_l\in c_0(S)$ and $\phi\in c_0(S)$, then for each $s, t\in S,$
\begin{eqnarray*}
\mid\frac{\phi\cdot\delta_s}{\omega}(t)\mid&=&\mid
\frac{\phi(st)}{\omega(st)}\mid\frac{\omega(st)}{\omega(t)}\\&\leq&
\parallel \phi\parallel_{\infty,\omega} w_l(t).
\end{eqnarray*}
 As $w_l\in c_0(S)$ then
$\frac{\phi\cdot\delta_s}{\omega}\in c_0(S)$, or equivalently. Similarly, from $w_r\in c_0(S)$ one may deduce that $\frac{\delta_s\cdot\phi}{\omega}\in
c_0(S).$ Therefore $\ell^1(S, \omega)$ is a dual Banach algebra, as required.

$(ii)$ The trivial  inequalities
 \[\frac{1_{t^{-1}s}}{\omega}\leq \frac{1}{\omega}~\emph{\emph{and}}~\frac{1_{st^{-1}}}{\omega}\leq \frac{1}{\omega},\ \ \ (s,t\in S), \]
 together with  $\frac{1}{\omega}\in c_0(S)$  imply that $\frac{1_{t^{-1}s}}{\omega},
\frac{1_{st^{-1}}}{\omega}\in c_0(S),$ and so
by Corollary~\ref{p3},   $\ell^1(S, \omega)$ is dual Banach
 algebra.
\end{proof}
As a partial converse to part $(ii)$ of the latter result we have,
\begin{proposition}\label{p5}
If  $\ell^1(S, \omega)$ is a  dual Banach algebra with predual $c_0(S)$ then  $\frac{1}{\omega}\in c_0(S)$, provided that either  $(st^{-1})^c$ or $(t^{-1}s)^c$ is finite, for some  $s,t\in S$.\\
 \end{proposition}
\begin{proof}
 Assume that $\ell^1(S, \omega)$ is dual Banach algebra with predual $c_0(S)$.  If  $(st^{-1})^c$ is finite for some $s, t\in S$, as $\frac{1_{st^{-1}}}{\omega}\in c_0(S)$ the function $\frac{1_{st^{-1}}}{\omega}$ vanishes on the complement of a finite subset $F$ of $S$. Therefore $\frac{1}{\omega}$ vanishes on  $((st^{-1})^c\cup F)^c$, which means that $\frac{1}{\omega}\in c_0(S)$.
\end{proof}
The next proposition describes a situation for which  $S$ is not weakly
cancellative and $\frac{1}{\omega}\notin c_0(S)$, nevertheless
$\ell^1(S, \omega)$ is dual Banach algebra with predual $c_0(S).$
\begin{proposition}\label{p6} Let $S_1$ be a weakly cancellative semigroup, $S_2$ be a non-weakly cancellative semigroup,
$\omega_1$ be a weight on $S_1$ such that $\frac{1}{\omega_1}$ is
not bounded and let $\omega_2$ be a weight on $S_2$ such that
$\frac{1}{\omega_2}\in c_0(S_2)$. Then $\ell^1(S, \omega_1\times\omega_2)$ is a dual
Banach algebra with respect to predual $c_0(S_1\times S_2)$, but
neither $S_1\times S_2$ is weakly cancellative  nor
$\frac{1}{ \omega_1\times\omega_2}\in c_0(S_1\times S_2).$
\end{proposition}
\begin{proof}
The identity \[\frac{1_{(s, t)(s', t')^{-1}}}{ \omega_1\times\omega_2}(x, y)=\frac{1_{ss'^{-1}}}{\omega_1}(x)\frac{1_{tt'^{-1}}}{\omega_2}(y)\ \ (x, y\in S),\] implies that $\ell^1(S, \omega_1\times\omega_2)$ is a left dual Banach algebra with predual $c_0(S_1\times S_2)$; indeed, for a given $\epsilon >0$,  as $S_1$ is weakly cancellative and $\frac{1}{\omega_2}\in c_0(S_2)$ there exist finite subsets $F_1\subseteq S_1$ and $F_2\subseteq S_2$ such that $\frac{1_{ss'^{-1}}}{\omega_1}(x)=0$ for $x\notin F_1$ and $\frac{1}{\omega_2}(y)<\epsilon$ for $y\notin F_2.$ Therefore, if $(x, y)\notin F_1\times F_2$ then $\frac{1_{(s, t)(s', t')^{-1}}}{\omega}(x, y)<M\epsilon$ where $M=\max_{x\in F_1}\frac{1_{ss'^{-1}}}{\omega_1}(x).$ A similar proof may apply for the right case.
That $\frac{1}{ \omega_1\times\omega_2}\notin c_0(S_1\times S_2)$ follows directly from the fact that $\frac{1}{\omega_1}$ is
not bounded.
\end{proof}

Recall that two weights $\omega$ and $\omega^{'}$ are said to be equivalent if there exist two positive numbers $\alpha, \beta$ such that $\beta\omega^{'}\leq\omega\leq\alpha\omega^{'},$ on the whole of $S$. If so then trivially  $\ell^1(S, \omega)=\ell^1(S, \omega^{'}).$   Two weights $\omega$ and $\omega^{'}$ on $S$ are said to be locally equivalent if for every $s, t\in S,$ there exist positive numbers $\alpha, \beta$ such that $\beta\omega^{'}\leq\omega\leq\alpha\omega^{'},$ on $st^{-1}\cup t^{-1}s$. Trivially equivalent weights are locally equivalent and the converse is also  hold, provided that  either $st^{-1}=S$ or $t^{-1}s=S$, for some $s,t\in S$. The next result compares  $\ell^1(S, \omega)$ and  $\ell^1(S, \omega^{'})$ from the dual Banach algebra point of view.
\begin{theorem}\label{p7} For every two locally equivalent weights $\omega$ and $\omega^{'}$ on $S$,  $\ell^1(S, \omega)$ is a dual Banach algebra with prdual $c_0(S)$ if and only if $\ell^1(S, \omega^{'})$ is  a dual Banach algebra with the predual $c_0(S).$
\end{theorem}
\begin{proof}Fix $s, t\in S$. As  $\omega$ and $\omega^{'}$ are locally equivalent, there exists positive numbers $\alpha, \beta$ such that $\beta\omega^{'}\leq\omega\leq\alpha\omega^{'}$ on $st^{-1}\cup t^{-1}s$. Let $\ell^1(S, \omega)$ be a dual Banach algebra with prdual $c_0(S).$ Then the equality  \[\frac{1_{st^{-1}}}{\omega^{'}}=\frac{1_{st^{-1}}}{\omega}\frac{\omega}{\omega^{'}}\]
together with the fact that $\frac{\omega}{\omega^{'}}\leq\alpha \ {\rm on }\ st^{-1}$ imply that $\frac{1_{st^{-1}}}{\omega^{'}}\in c_0(S);$ see Corollary~\ref{p3}. Similarly  $\frac{1_{t^{-1}s}}{\omega^{'}}\in c_0(S)$ and so $\ell^1(S, \omega^{'})$ is  a dual Banach algebra with prdual $c_0(S).$

\end{proof}

We conclude this section  with the following examples illustrating the last results.
\begin{example}\label{1}

$(i)$ Let $S=(\mathbb{N}, \vee)$, where $m\vee n=\max\{m,n\}$. As $S$ is weakly cancellative,
$\ell^1(S,\omega)$ is a dual Banach algebra with  predual
$c_0(S)$.

$(ii)$  Let $S=(\mathbb{N},\wedge),$ where $m\wedge n=\min\{m,n\}$. Then $1\wedge 1^{-1}={\mathbb{N}}$ and so  $S$ is not weakly cancellative. Since  $(1\wedge 1^{-1})^c$ is empty, Proposition~\ref{p5} together with part $(ii)$ of Proposition~\ref{p4} imply that $\ell^1(S, \omega)$ is a dual Banach algebra with predual $c_0(S)$ if and only if  $\frac{1}{\omega}\in c_0(S)$. Note that two weights $\omega(n)=1+n$ and $\omega^{'}(n)=1+n^2$ are not locally equivalent on $S$; nevertheless both of  $\ell^1(S, \omega)$ and  $\ell^1(S, \omega^{'})$ are dual Banach algebra with predul $c_0(S)$, and this confirms that the converse of Theorem~\ref{p7} is not valid, in general.
\end{example}
 \begin{example}\label{2}
 Let $S=(\mathbb{Z},\wedge)$, $\omega$ be a weight on $S$ and let $\omega_{\mathbb N}$ be the restriction of $\omega$ on the subsemigroup $\mathbb{N}$ of $S$. A direct verification shows that for $m, n\in S$,
  \begin{eqnarray*} m\wedge n^{-1}&=&\left\{%
\begin{array}{lll}
    \phi & \hbox{$m>n$} \\
    \{m\} & \hbox{$m<n$} \\
\{m, m+1, \cdots \} & \hbox{$m=n$}. \\
\end{array}%
\right.
\end{eqnarray*}
\noindent In particular, $1\wedge 1^{-1}=\mathbb{N}$. If $\ell^1(S, \omega)$ is dual Banach algebra with predual $c_0(S)$ then, by Corollary~\ref{p3}, $\frac{1_\mathbb{N}}{\omega}=\frac{1_{(1\wedge 1^{-1})}}{\omega}\in c_0(S)$ from which we have $\frac{1}{\omega_{\mathbb N}}\in c_0({\mathbb N})$. Conversely, if  $\frac{1}{\omega_{\mathbb N}}\in c_0({\mathbb N})$ then for each $\epsilon>0$ there exists  $N\in\mathbb N$ such that $\frac{1}{\omega_{\mathbb N}}(k)<\epsilon$, for all $k>N$. This implies that for each $m, n\in S$, if $k\in\mathbb Z$ is so that $\mid k\mid>\max\{N, n\}$ then $\frac{1_{(m\wedge n^{-1})}}{\omega}(k)<\epsilon;$ or equivalently, $\frac{1_{(m\wedge n^{-1})}}{\omega}\in c_0(S)$. Therefore  Corollary~\ref{p3} implies that $\ell^1(S, \omega)$ is dual Banach algebra with predual $c_0(S).$ We have thus proved that for every weight $\omega$ on $S=(\mathbb{Z},\wedge)$ the next assertions are equivalent:

 $(i)$  $\ell^1((\mathbb{Z},\wedge), \omega)$ is dual Banach algebra with predual $c_0(\mathbb Z);$

 $(ii)$  $\ell^1((\mathbb{Z},\vee), \widetilde{\omega})$ is dual Banach algebra with predual $c_0(\mathbb Z)$, where $\widetilde{\omega}(n)=\omega(-n);$

$(iii)$  $\ell^1(\mathbb N, \omega_{\mathbb N})$ is dual Banach algebra with predual $c_0(\mathbb N);$

$(iv)$  $\frac{1}{\omega_{\mathbb N}}\in c_0({\mathbb N}).$\\
In which, the equivalence $(i)\Leftrightarrow(ii)$ follows from the fact that the Banach algebras  $\ell^1((\mathbb{Z},\wedge), \omega)$  and  $\ell^1((\mathbb{Z},\vee), \widetilde{\omega})$ are isometrically isomorphic and equivalence $(iii)\Leftrightarrow(iv)$ is proved in part $(ii)$ of Example~\ref{1}.

One may trivially accept that $(iv)$ dose not imply   that $\frac{1}{\omega}\in c_0(\mathbb Z)$, in general; as a sample, for the weight \begin{eqnarray*}\omega(n)&=&\left\{%
\begin{array}{ll}
    1+n & \hbox{$n\in\mathbb N$} \\
    1 & \hbox{${\rm otherwise}$}, \\
\end{array}%
\right.
\end{eqnarray*}
$\ell^1((\mathbb{Z},\wedge), \omega)$ is dual Banach algebra with predual $c_0(\mathbb Z)$, but  $\frac{1}{\omega}$ does not belong to $c_0(\mathbb Z)$. Furthermore  $({m\wedge n^{-1}})^c$ and $({n\wedge m^{-1}})^c$ are infinite for all $m, n\in\mathbb Z$, and this confirms that the hypothesis imposed  in Proposition~\ref{p5} are essential.
\end{example}

\begin{example}\label{3}
 Let $S=({\mathbb{N}}\times {\mathbb{N}},\cdot)$, where
$(m, n)\cdot(m',n')=(m+m',n')$ and let $\omega(m,n)=e^{-m}(1+n)$, for
all $m, n, m', n'\in {\mathbb{N}}$.  Then  $\ell^1(S, \omega)$
is a dual Banach algebra with respect to predual $c_0(S)$. However,
$S$ is not weakly cancellative (for instance, $(2, 2)\cdot(1, 2)^{-1}=1\times\mathbb{N}$), and $\frac{1}{\omega}\notin
c_0(S).$ This is a modification of Proposition~\ref{p6} for $\omega_1(m)=e^{-m}$, $\omega_2(n)=1+n$, $S_1= (\mathbb{N}, +)$ and  $S_2=\mathbb{N}$, equipped with the right zero multiplication.
\end{example}

\begin{example}\label{4} Let $S$ be an infinite left zero semigroup,  then trivially  $s^{-1}s=S$ and $st^{-1}=\{s\}$, for each $s, t\in S$. Therefore $S$ is not left weakly cancellative however it is right weakly cancellative; and so  for every weight $w$ on $S$,  $\ell^1(S, \omega)$ is a left dual Banach algebra with predual $c_0(S)$, but it is not a right dual Banach algebra. It is worthwhile mentioning that  $\ell^1(S)$ is actually a dual Banach algebra with respect to the predual $c(S)=c_0(S)\oplus{\mathbb C}1,$ (see {\cite[Example 7.30]{Dales-Lau-s}}).
\end{example}

\section{semigroup measure algebra as a dual Banach algebra}
Throughout this section  $S$ is assumed to be, at least, a locally compact Hausdorff topological  semigroup. The (semigroup)
measure algebra $M_b(S)$ of $S$ is the space of all bounded Radon measures
on $S$. It is known that $M_b(S)$ is the dual of $C_0(S)$ as Banach space (under the duality $\langle \mu, f\rangle=\int_Sfd\mu$, for $\mu\in M_b(S), f\in C_0(S)$), where $C_0(S)$ is the space of all bounded continuous functions on $S$ which vanishes  at infinity.  Furthermore, $M_b(S)$ is  a Banach algebra under the following convolution product.
$$\mu *\nu(f)=\int_S\int_Sf(st)d\mu(s)d\nu(t),\ \ \ \ \ (f\in C_0(S), \mu,\nu\in M_b(S)).$$
A direct verification reveals that the induced module actions are
$$(f\cdot\mu)(s)=\int_Sf(st)d\mu(t)\ \ \ {\rm and} \ \ (\mu\cdot f)(s)=\int_Sf(ts)d\mu(t), \ \ \ \ \ \ \ (\mu\in M_b(S),  f\in C_0(S), s\in S).$$

The next result impose some conditions on $S$ under which $M_b(S)$ becomes a dual Banach algebra with predual $C_0(S)$, (see also  {\cite[Proposition 2.3]{pym}}).
\begin{proposition}\label{M} If the sets  $FK^{-1}$ and $K^{-1}F$ are relatively compact, for all compact subsets $F$ and $K$ of $S$, then the measure algebra $M_b(S)$
is a dual Banach algebra with respect to predual $C_0(S)$.\\
\end{proposition}
\begin{proof}
 Let  $f\in C_0(S)$,  $\mu\in
M_b(S)$ and $\epsilon>0$  be arbitrary. Therefore there exist  compact
subsets $F$ and $K$ of  $S$ such that
$$|f(s)|<\epsilon ~~\emph{\emph{for all}}~~ s\notin F \ \ {\rm and}\ \ |\mu(S\backslash K)|<\epsilon.$$ Then for  $s\notin\overline{K^{-1}F}, $ (which is compact by the hypothesis)
\begin{eqnarray*}|\mu\cdot f(s)|&=&|\int_Sf(ts)d\mu(t)|\\&\leq&|\int_Kf(ts)d\mu(t)|+|\int_{S\backslash K}f(ts)d\mu(t)|\\&<&(\|\mu\|+\|f\|_\infty)\epsilon.\end{eqnarray*}
 Therefore $\mu\cdot
f\in C_0(S)$. A similar argument may apply for  $f\cdot\mu\in C_0(S).$
\end{proof}
  Since all compact Hausdorff spaces and also all locally compact Hausdorff topological groups  trivially obey the hypothesis of the latter theorem,  as an immediate consequence  we have the next corollary.
\begin{corollary}\label{M1} If $S$ is either compact  or a  topological group, then  $M_b(S)$ is a dual
Banach algebra with  predual $C_0(S).$
\end{corollary}
 The next result provides a partial converse to Theorem~\ref{M}.
 \begin{proposition}\label{N} If $M_b(S)$ is a dual Banach algebra with predual $C_0(S)$ then $s^{-1}K$ and $Ks^{-1}$ are  compact,  for each compact set $K$ in $S$ and $s\in S.$
 \end{proposition}
\begin{proof}
Let $K$ be a compact subset of $S$ and let $s\in S$, then  there exists an $f\in C_0(S)$ such that $f(K)=1.$ Since $M_b(S)$ is a dual Banach algebra with predual $C_0(S)$,  $\delta_s\cdot f\in C_0(S)$ and thus $\{t\in S: f(st)=(\delta_s\cdot f)(t)\geq 1\}$ is compact and contains $s^{-1}K
$ as a closed subset. Therefore $s^{-1}K$ is  compact. Similarly, $Ks^{-1}$ is  compact.
\end{proof}
It should be remarked that, the compactness of  $s^{-1}K$ and $Ks^{-1}$,  for every compact set $K$ in $S$ and $s\in S$, is also equivalent to the fact that $C_0(S)$ is translation invariant; for more details see {\cite[Example 3.1.10]{BJM}}, in which one may find some other equivalent assertions of this fact from the semigroup compactification point of veiw.\\

Using Corollary~\ref{M1} and Proposition~\ref{N} we have the next result.
\begin{corollary}
Let $S$ be either a left zero, a right zero or a zero semigroup, then $M_b(S)$ is a dual Banach algebra with predual $C_0(S)$ if and only if $S$ is compact.
\end{corollary}
\begin{proof}
Follows trivially from the above mentioned results and the fact that, in the  right (respectively, left) zero semigroup setting, $ss^{-1}=S$ (respectively, $s^{-1}s=S$), for each $s\in S$; and in the case where $S$ is a zero semigroup (with a zero $z$) then $zs^{-1}=S=s^{-1}z$, for each $s\in S$.
\end{proof}

\section*{acknowledgement}
 The second author is partially supported by the Center of Excellence for Mathematics at Isfahan
university and the third one is partially supported by a research grant from  Ferdowsi University of
Mashhad.


\begin{thebibliography}{}
\bibitem{BJM}{ J. F. Berglund, H. D. Junghenn \and P. Milnes,} {\em
Analysis on Semigroups}, Wiley, New York, (1989).
\bibitem{CY} {P. Civin \and B. Yood}, {\em The second conjugate space of a Banach algebra as an algebra}, Pacific. J. Math.,
\textbf{11} (1961), 847--870.
\bibitem{DL}{H. G. Dales \and A. T.-M. Lau}, {\em The second duals of Beurling algebras}, Mem. Amer. Math. Soc. \textbf{177} (2005), no. 836.
\bibitem{Dales-Lau-s} H. G. Dales, A. T.-M. Lau \and D. Strauss, {\em Banach algebras on semigroup and their compactifications}, Preprint.
\bibitem{DPV}
{ H. G. Dales, A. Rodrigues-Palacios \and M. V. Velasco}, {\em The second transpose of a derivation},
{ J. London Math. Soc.},  \textbf{64}(2) (2001), 707--721.
\bibitem{Daws} M. Daws, {\em Connes-amenability of bidual and weighted semigroup
algebras}, Math. Scand. \textbf{99} (2006) no. 2, 217-246.
\bibitem{DPW} M. Daws, H. L. Pham \and S. White, {\em Conditions implying the uniqeness of the weak$^*-$topology on certain group algebras}, arXiv:[math.FA]0804.3764v1.
\bibitem{MV} S. Mohammadzadeh \and H. R. E. Vishki, {\em Arens regularity of  module actions and the second adjoint  of a derivation.} To appear in Bull. Austral. Math. Soc. (2008).
\bibitem{pym} J. S. Pym, {\em Idempotent measures on semigroups}, Pacific J. Math. \textbf{12} (1962), 685-698.

\bibitem{Runde1} V. Runde, {\em Lectures on amenability}, Lecture Notes in Mathematics, \textbf{1774}
Springer-Verlag, Berlin, (2002).
\bibitem{Rundea} V. Runde, {\em Amenability for dual Banach algebras},
Studia Math. \textbf{148}, (2001), 47-66.
\bibitem{Y} {N. J. Young}, {\em The irregularity of multiplication in group algebras},  Quart. J. Math. Oxford \textbf{24} (2),  (1973), 59-62.
\end{thebibliography}
\end{document}